\documentclass{ifacconf}

\newtheorem{nnassumption}{\bf Assumption}

\newtheorem{nntheorem}{\bf Theorem}
\newenvironment{theorem}{\begin{nntheorem}\it}{\end{nntheorem}}
\newtheorem{nndefinition}{\bf Definition}

\newtheorem{nnproposition}{\bf Proposition}

\newtheorem{nnproblem}{\bf Problem}

\newtheorem{nnlemma}{\bf Lemma}
\newenvironment{lemma}{\begin{nnlemma}\it}{\end{nnlemma}}
\newtheorem{nnremark}{\bf Remark}
\newenvironment{remark}{\begin{nnremark}}{\hfill \hspace*{1pt}\hfill $\circ$\end{nnremark}}
\newenvironment{proof}{{\bf Proof.}}{\hfill \hspace*{1pt}\hfill $\Box$}

\usepackage{amsfonts,latexsym,amsmath,amssymb}
\usepackage[mathscr]{eucal}
\usepackage{graphicx,color}
\usepackage{comment}

\newcommand\eps{\varepsilon}

\usepackage{graphicx}      
\usepackage{natbib}        
\begin{document}
\begin{frontmatter}

\title{Event-triggered damping stabilization \\of a linear wave equation\thanksref{footnoteinfo}} 

\thanks[footnoteinfo]{This work was partly funded by the ERC Advanced Grant Taming.}

\author[First]{Lucie Baudouin} 
\author[First]{Swann Marx} 
\author[First]{Sophie Tarbouriech}

\address[First]{LAAS-CNRS; Universit\'e de Toulouse; CNRS, Toulouse; France.}

\begin{abstract}                
The paper addresses the design of an event-triggering mechanism for a partial differential wave equation posed in a bounded domain. The wave equation is supposed to be controlled through a first order time derivative term distributed in the whole domain. Sufficient conditions based on the use of suitable Lyapunov functional are proposed to guarantee that an event-triggered distributed control still ensures the exponential stability of the closed-loop system. Moreover, the designed event-triggering mechanism allows to avoid the Zeno behavior. The `existence and regularity' prerequisite properties of solutions for the closed loop system are also proven.
\end{abstract}


\end{frontmatter}

\section{Introduction}
In this paper we are concerned with a stabilization problem by event-based control of the wave equation 
\begin{equation}
	\partial_t^2 z(t,x) - \Delta z (t,x) = u(t,x),\quad x \in \Omega, \quad t\geq 0\\
\end{equation}
in a suitable functional space, where $\Omega$ is a bounded domain in $\mathbb R^N$, $\Delta$ is the Laplacian operator, the Dirichlet boundary datum of the unknown $z$ is homogeneous and $u=u(x,t)$ is a distributed parameter control that can be sampled in the time variable. 
One of the most classical way to obtain stabilization and even exact controllability is to choose $u(x,t) = - \alpha \partial_t z(x,t)$ that corresponds to a classical damping term (see \cite{ChenSICON79}).

There exist many ways to stabilize or to analyze the stability of the wave equation. For the analysis of the stability, let us cite the seminal paper \cite{ChenSICON79}, which considers a full damping and uses some Lyapunov technique to characterize the stability, or \cite{lebeau1996equation}, which focuses on a localized damping and follows a micro-local analysis to provide the optimal decay rate of the solution. The multiplier method developed in \cite{komornik1997exact} is also a powerful tool to prove controllability and stabilization of several partial differential equations (PDEs). For the feedback design, we might also mention \cite{smyshlyaev2010boundary} in which the backstepping method is used for several kinds of one dimensional wave equations.

In the context of event-triggered control, two objectives can be pursued: (1) Emulation - the controller is a priori predesigned and only the event-triggered rules have to be designed (see, for example, \cite{pos:tab:nes:ant/IEEE2015}, \cite{esp:tan:tar/CDC2017} and the references therein), or (2) Co-design - the joint design of the control law and the event-triggering conditions has to be performed (see, for example, \cite{seuret2016lq}, \cite{heemels2013periodic}  and the references therein). 
The current paper comes within the scope of the first case, that is, of the emulation context and, more precisely, we want to focus on the stabilization of the wave equation by means of an event-triggered control. 
This kind of feedback strategy has been followed in order to reduce the computational resources of the controllers. Indeed, event-triggered controls can be defined as controls updated aperiodically, only when it is needed, i.e when the system becomes, roughly speaking, too unstable. 
This strategy has been first applied for finite-dimensional systems (see \cite{aastrom1999comparison} \cite{tabuada2007event} or \cite{peralez2018event}) and then extended to some infinite-dimensional systems (see \cite{espitia2016event} for the case of boundary control of linear hyperbolic system of first order PDEs and \cite{selivanov2016distributed} for the case of distributed control of some diffusion equations). 
Let us also mention \cite{davo2018sample}, where a first order hyperbolic linear PDE is considered with a sampled-data control, which is a control strategy closely related to the event-triggered control. 
To the best of our knowledge, the current paper is the first one dealing with a wave equation with an event-triggered distributed control. 

Let us emphasize that in the event-triggered control context, the plant evolves in continuous time, whereas the control signal is updated depending on discrete-time events. 
The resulting closed-loop system hence becomes a hybrid system in the sense given in \cite{teel/book2012}. Note that coupling these two components requires to be careful. Indeed, such a controller might produce a Zeno solution, which is defined as ``jumping" infinitely in a bounded time interval - read \cite[Chapter 2]{teel/book2012} for a complete discussion on this topic. This kind of solutions has to be avoided, otherwise the implementation of the feedback law might be impossible. 
Here, instead of considering the hybrid framework to develop the results, we follow an alternative route as in \cite{tabuada2007event} and 
\cite{esp:tan:tar/CDC2017}.

In the current paper, the design of an event-triggering mechanism for a wave equation posed in a bounded domain is proposed. 
We consider the case of a linear wave equation submitted to an event-triggered distributed control. In particular, under some assumptions linked with  the shape and size of the spatial domain $\Omega$, the control mechanism is designed to avoid Zeno solution and ensure the exponential stability of the closed-loop system. The originality of the approach relies on the conditions involving the parameters of the event-triggering mechanism and the characterization of the overshoot estimation and the  exponential rate (defining the exponential stability property). The  results are developed thanks to the use of an adequate Lyapunov functional.

This paper is organized as follows. Section~\ref{sec_state} describes the system under consideration and states the problem to be solved  in a more precise way. Then,  the main results of the paper in terms of well-posedness and exponential stability theorems are presented in Section~\ref{sec_main}. Section~\ref{sec_wp} is devoted to the proof of the well-posedness of the closed-loop system and Section~\ref{sec_exp} focuses on the proof of the exponential stability using Lyapunov theory. Finally, Section~\ref{sec_discussion} collects some concluding remarks and gives further research lines to be followed in the future.

{\bf Notation: } When applicable, the partial derivative with respect to $t$ is denoted $\partial_t z = \frac{\partial z}{\partial t}$, the gradient is the vector function $\nabla z = (\partial_{x_n} z)_{n=1,...,N}$ and the Laplacian operator is defined by  $\Delta z =  \sum_{n=1}^{N}\partial_{x_n}^2 z $.
The Hilbert space of square integrable functions over $\Omega$ with values in $\mathbb R$ is denoted by $L^2(\Omega)$ with norm $\|z\|_{L^2(\Omega)} = \left(\int_\Omega |z(x)|^2 dx\right)^{1/2} $ and the Sobolev space $H^1(\Omega)$ is  $\big\{ z\in L^2(\Omega),  \nabla z \in \left(L^2(\Omega)\right)^N \big\}$ with $\|z\|^2_{H^1(\Omega)} = \|z\|_{L^2(\Omega)}^2 + \|\nabla z\|_{L^2(\Omega)}^2$.
We denote finally $C_0(\mathbb R_+)$ the space of continuous functions defined over $\mathbb R_+$ and vanishing at infinity.

\section{Problem statement}
\label{sec_state}
  
The goal of our work is to study the stabilization of the partial differential system
\begin{equation}\label{EqU}
	\left\{
	\begin{aligned}
		&\partial_t^2z(t,x) - \Delta z(t,x) = u(t,x),\qquad (t,x)\in\mathbb{R}_+\times \Omega\\
		&z(t,x)=0,\qquad \qquad \qquad \qquad\quad  ~  (t,x)\in\mathbb{R}_+\times \partial \Omega\\
		&z(0,x)=z_0(x),\: \partial_t z(0,x)=z_1(x),\qquad  \qquad x\in \Omega,
	\end{aligned}
	\right.
\end{equation}
where $z$ and $u$ denote the state and the control variable respectively. 
Before going any further, let us recall that if $u \in L^1(\mathbb R_+ ;L^2(\Omega))$, 
 $(z_0,z_1)\in H^1(\Omega)\times L^2(\Omega)$ and assuming the compatibility condition $z_0(x) = 0$ for all $x\in\partial\Omega$ (so that $z_0$ actually belongs to the Hilbert space denoted $H^1_0(\Omega)$), then the initial and boundary value problem \eqref{EqU} is well-posed and has a unique solution (see e.g. \cite{Lions}, \cite{LasieckaLionsTriggiani}, \cite{Haraux})
$$
	 z \in C^0(\mathbb R_+; H^1(\Omega)) \cap C^1(\mathbb R_+; L^2(\Omega)).
$$ 

If one sets $u(t,x):=-\alpha \partial_t z(t,x)$, where $\alpha >0$, then the origin of the closed-loop system is globally exponentially stable (see e.g. \cite{ChenSICON79}, \cite{ChenSICON81}).
We can build an explicit Lyapunov functional for such a system, which is closely related to the natural energy (sum of the kinetic and potential energies) of the wave equation given by
\begin{equation}\label{energy}
	\begin{aligned}
		E(t)&:=\frac{1}{2} \left( \|\partial_t z(t)\|_{L^2(\Omega)}^2 + \|\nabla z(t)\|_{L^2(\Omega)}^2 \right)
	\end{aligned}
\end{equation}
In this paper, the objective relies on the way to implement the control input $u=-\alpha \partial_t z$, taking into account that $u$ is performed through a sample-and-hold mechanism. It is not continuously updated or transmitted to the actuators. Indeed, it is updated at certain instants $\{t_k\}_{k\in\mathbb N}$, which form a sequence of strictly increasing positive scalars. Control action is held constant between two successive sampling instants ($t_k$ and $t_{k+1}$) through a zero order holder. Differently from classical digital control approaches, the sampling interval $t_{k+1}-t_k$ is not assumed to be
constant. Indeed, the control input can be written as 
\begin{equation}
	u(t,x) = -\alpha \partial_t z(t_k,x),\quad \forall (t,x)\in [t_k,t_{k+1}[\times \Omega .
\end{equation}

The idea is to find an appropriate event-triggered law for the sampling times $(t_k)_{k\in\mathbb{N}}$ in order to maintain the exponential stabilization toward the origin of the closed loop 
\begin{equation}\label{Eqet}
	\left\{
	\begin{aligned}
		&\partial_t^2z(t,x) - \Delta z(t,x) = -\alpha \partial_t z(t_k,x),\\
		&\quad \qquad \qquad ~ \qquad \qquad (t,x)\in[t_k,t_{k+1}[\times \Omega, \forall k\in\mathbb N\\
		&z(t,x)=0,\qquad \qquad \qquad \qquad \quad (t,x)\in\mathbb{R}_+\times \partial \Omega\\
		&z(0,x)=z_0(x),\: \partial_t z(0,x)=z_1(x),\qquad \qquad x\in \Omega.
	\end{aligned}
	\right.
\end{equation}

Besides, since the sampling will be aperiodic, we will have to carefully check whether such a feedback law does not produce any Zeno phenomenon.
Denoting by $e_k$ the speed deviation, given by the formula
\begin{equation}\label{deviation}
	e_k(t,\cdot):=\partial_t z(t,\cdot)-\partial_t z(t_k,\cdot),~ \hbox{ in } \Omega, ~ \forall t\in [t_k,t_{k+1}[,
\end{equation}
the sampling instants are determined from the following logic:
\begin{equation}\label{et-law}
	\begin{aligned}
		t_{k+1}:=\inf\Big\{ t\geq t_k, &~\|e_k(t)\|_{L^2(\Omega)}^2 - \gamma_0  \|z(t)\|_{L^2(\Omega)}^2 \\
		&~-\gamma_1 \|\partial_t z(t)\|_{L^2(\Omega)}^2  - \eta_0(t) \geq 0\Big\},
	\end{aligned}
\end{equation}
where $\gamma_0$ and $\gamma_1$ are sufficiently small positive constant that have to be defined and $\eta_0$ is a function that will be specified later, decreasing and strictly positive. 

\begin{remark}
We could have designed the event-triggered law more simply with $\gamma  E(t)$
but we wanted to keep a close track on the respective constraints we need to impose in position (on $\gamma_0$) and velocity ($\gamma_1$) of the variable $z$.
\end{remark}

\section{Main Results}
\label{sec_main}

This section collects the two main theorems of the paper. The first one states the existence and classical regularity of solutions for the closed-loop event-triggered controlled system \eqref{Eqet}-\eqref{et-law} and addresses the Zeno behavior.

\begin{theorem}[Well-posedness]\label{WP}
\label{thm-wp}
Consider the linear wave equation \eqref{Eqet} under the event-triggering mechanism \eqref{et-law}. 
For any initial condition $(z_0,z_1)\in H^1_0(\Omega)\times L^2(\Omega)$, there exists a unique solution 
	$$ z \in C^0(\mathbb R_+; H^1(\Omega)) \cap C^1(\mathbb R_+; L^2(\Omega)). $$
Furthermore, the Zeno phenomenon is avoided. 
\end{theorem}


The second result states the global exponential stability of the origin of the closed-loop system \eqref{Eqet}-\eqref{et-law} under a specific condition on the spatial domain. 
\begin{theorem}[Global exponential stability]\label{GES}
Let $\alpha>0$ be the damping parameter in the wave equation closed-loop system \eqref{Eqet}. 
Assume that the Poincar\'e constant of the bounded domain $\Omega$, denoted $C_\Omega$ (see Lemma~\ref{Poincare}), satisfies 
$$C_\Omega < \sqrt 2.$$
There exist $\gamma_0>0$, $\gamma_1>0$ and a strictly decreasing function $\eta_0 \in C_0(\mathbb R_+)$ such that for any initial condition $(z_0,z_1)\in H^1_0(\Omega)\times L^2(\Omega)$, system \eqref{Eqet} under the event-triggering mechanism \eqref{et-law} is exponentially stable:
  \begin{equation}\label{ExpStab}
			\exists K, \delta >0 \hbox{ such that }  E(t) 	\leq K e^{-\delta t} E(0) 
			, \quad  \forall t>0.
	\end{equation}
\end{theorem}
\begin{remark}
Note that the constants $\gamma_0$ and $\gamma_1$, and the function $\eta_0$ can be built explicitely. This construction is given in the fourth step  of the proof of the global exponential stability, that is written in Section \ref{sec_exp}. 
\end{remark}
\begin{remark}
\label{rem-restriction-domain}
The restriction on the Poincar\'e constant $C_\Omega$ (see Lemma~\ref{Poincare}) given in Theorem \ref{GES} is closely related to the shape of the domain $\Omega$. It is important to note that the condition can actually be viewed as a link between the domain $\Omega$ and the velocity of the wave equation under consideration. For simplicity, we have chosen in this paper a wave velocity identically equal to $1$. In a more general case, the wave equation may be rewritten as 
$$
\partial_t^2 z(t,x) - c^2 \Delta z(t,x) = u(t,x),\quad \forall (t,x)\in \mathbb{R}_+ \times \Omega,
$$
where $c$ denotes the velocity, and when studying this equation, the restriction on the Poincar\'e constant  becomes $c^{-1}C_\Omega < \sqrt 2$. Noticing that $C_\Omega$ is closely related to the diameter of the domain $\Omega$, the latter comment implies, roughly speaking, that a too small velocity $c$ will not imply a stabilization result with such an event-triggered damping control if the domain is too large. In other words, our stabilization result is limited to ``small" ratio between the size of the domain and the velocity of the wave equation (and not only to ``small" constants $C_\Omega$).    
\end{remark}

\section{Well-posedness of the closed loop system}
\label{sec_wp}
The proof of Theorem~\ref{WP} is divided into three steps. We first prove that the closed-loop system \eqref{Eqet}-\eqref{et-law} is well-posed on every sample interval $[t_k,t_{k+1}]$ in a way such that one obtains a unique solution $z\in C([0,T],H^1_0(\Omega))\cap C^1([0,T];L^2(\Omega))$ for any $T>0$. 
Then, we show that \eqref{Eqet}-\eqref{et-law} avoids the Zeno phenomenon. Finally, gathering these information proves that the solution $z$ evolves in $\mathbb R_+\times\Omega$. 

\begin{proof} 

\noindent $\bullet$ {\it Existence, uniqueness and regularity of the solution:} \\
We proceed by induction. First, let us focus on the initialization interval $[0,t_1]$, and prove that the solution belongs to the awaited functional space and is unique. Then, we will assume that for a fixed integer $k$ the regularity holds true up to $t_{k+1}$, and proceed on the next time interval. 

\noindent $(i)$ {\it Initialization.} On the first time interval, \eqref{Eqet} reads
$$
\left\{
\begin{array}{lr}
\partial_t^2 z(t,x) - \Delta z(t,x) = -\alpha z_1(x),& (t,x)\in [0,t_1]\times \Omega,\\
 z(t,x) = 0,&  (t,x)\in [0,t_1]\times \partial\Omega,\\
 z(0,x) = z_0(x),\: \partial_t z(0,x) = z_1(x),& x\in \Omega.
\end{array}
\right.
$$
This is a simple wave equation with source term. By assumption, $z_1 \in L^2(\Omega)$ so that $-\alpha z_1 \in L^1(0,t_1;L^2(\Omega))$, and one can invoke \cite[Theorem 1.3.2]{Haraux} to conclude that there exists a unique solution $z\in C([0,t_1];H^1_0(\Omega))\cap C^1([0,t_1];L^2(\Omega))$ to the latter system. 

\noindent $(ii)$ {\it Heredity.} Fix $k\in \mathbb{N}$ and assume that $$z \in C([t_k,t_{k+1}];H^1_0(\Omega))\cap C^1([t_k,t_{k+1}];L^2(\Omega)).$$ 
Denote by $z_{2k+2}$ and  $z_{2k+3}$ the position and velocity function values of the wave at time $t_{k+1}$.
Now consider \eqref{Eqet} over the time interval $[t_{k+1},t_{k+2}]$:
$$
\left\{
\begin{aligned}
&\partial_t^2 z - \Delta z =  - \alpha z_{2k+3} ,\qquad\qquad~ \hbox{in } [t_{k+1},t_{k+2}]\times \Omega,\\
& z = 0,\qquad\qquad\qquad\qquad \qquad  ~\hbox{on } [t_{k+1},t_{k+2}]\times \partial\Omega,\\
&z(t_{k+1}) = z_{2k+2},\: \partial_t z(t_{k+1}) = z_{2k+3},\qquad ~ \qquad\hbox{in } \Omega.
\end{aligned}
\right.
$$
It is again a wave equation with source term, here $- \alpha z_{2k+3}$, that belongs to $L^1([t_{k+1},t_{k+2}];L^2(\Omega))$ since we assumed 
$\partial_t z \in  C([t_{k},t_{k+1}];L^2(\Omega))$ and $\partial_t z(t_{k+1}) = z_{2k+3}$. Thus, applying \cite[Theorem 1.3.2]{Haraux} we conclude to the existence and uniqueness of the solution 
$$
z \in C^0([t_{k+1},t_{k+2}];H^1_0(\Omega))\cap C^1([t_{k+1},t_{k+2}];L^2(\Omega)).
$$

\noindent $(iii)$ {\it Conclusion.} By induction, for any $k\in\mathbb{N}$, 
$
z \in C([t_{k},t_{k+1}];H^1_0(\Omega))\cap C^1([t_{k},t_{k+1}];L^2(\Omega)).
$
Therefore, from the extension by continuity at the instants $t_k$, one can conclude that \eqref{Eqet} has a unique solution 
\begin{equation*}
z \in C([0,T];H^1_0(\Omega))\cap C^1([0,T];L^2(\Omega)).
\end{equation*}  

\noindent $\bullet$ {\it Avoiding Zeno phenomenon:} \\
The goal is to prove that with  the event-triggering mechanism \eqref{et-law}, there is only a finite number of sampling times over any closed interval subset
 $[0,T]$ of $\mathbb R_+$. 
 In fact, given $T>0$, we will show that there exists $\tau^*>0$ such that the sampling times $\{t_k\in [0,T], k\in\mathbb N\}$ from \eqref{et-law}  satisfy $t_{k+1}- t_k > \tau^*$. 
 
The proof relies on the continuity of $t\mapsto \partial_t z(t,\cdot)$ as a function from $[0,T]$ to $L^2(\Omega)$.
The uniform continuity of the 
function $t\mapsto e_k(t,\cdot)$ defined by \eqref{deviation} from the compact set $[0,T]$ to the Hilbert space $L^2(\Omega)$ stems from $\partial_t z \in C([0,T] ; L^2(\Omega))$. The contrapositive of the definition of this uniform continuity brings that 
$\forall \eta >0$, $\exists \tau^*>0$, $\forall s,t \in [0,T]$
$$
\|e_k(t) - e_k(s)\|_{L^2(\Omega)} > \eta
\quad \Rightarrow \quad
|t-s| > \tau^*.
$$

Hence, applying this property to $s = t_k$, $t= t_{k+1}$, we ultimately need to prove that we can bound from below, independently of $k$, the quantity 
$$\|e_k(t_{k+1}) - e_k(t_k)\|_{L^2(\Omega)} = \|e_k(t_{k+1})\|_{L^2(\Omega)}$$
 in order to avoid the Zeno phenomenon. Now, by definition of $t_{k+1}$ in \eqref{et-law} and $\eta_0$ (decreasing and strictly positive), we have
$$\begin{aligned}
	&\|e_k(t_{k+1})\|_{L^2(\Omega)}^2 \\
	&\geq \gamma_0  \|z(t_{k+1})\|_{L^2(\Omega)}^2 + \gamma_1 \|\partial_t z(t_{k+1})\|_{L^2(\Omega)}^2 
	+ ~\eta_0(t_{k+1}) \\
	&\geq \eta_0(T) > \eta > 0,
\end{aligned}$$
that allows to conclude that the Zeno behavior is avoided.

\noindent $\bullet$ {\it Conclusion:}\\
A unique solution exists on any time interval $[0,T]$, for all strictly positive $T$. Moreover, the Zeno phenomenon is avoided in such bounded intervals, so that in the end, we can write
$z \in C^0(\mathbb{R}_+;L^2(\Omega))\cap C^1(\mathbb{R}_+;L^2(\Omega))$, allowing to end the proof of Theorem~\ref{WP}.
\end{proof}

\section{Exponential stability}
\label{sec_exp}
This section is devoted to the proof of Theorem \ref{GES}. Let us consider the following Lyapunov functional candidate :  
\begin{multline}\label{Lyap}
V(t):= \frac{1}{2}\int_\Omega |\partial_t z(t,x)|^2 dx + \frac{1}{2}\int_\Omega |\nabla z(t,x)|^2 dx \\
+ \frac{\eps\alpha}{2}\int_\Omega |z(t,x)|^2 dx + \eps\int_\Omega  z(t,x)\partial_t z(t,x) dx,
\end{multline}
defined with $\eps>0$ and for the state variable $z(t)\in L^2(\Omega)$ of \eqref{Eqet} and where
we denoted the functional by $V(t)$ instead of $V(z(t))$ in sake of simplicity. Then, one gets: 
	$$V(t)= E(t) + \dfrac {\eps\alpha}2  \|z(t)\|_{L^2(\Omega)}^2+ \eps\int_\Omega  z(t)\partial_t z(t).$$
The proof is divided into four steps (for all $t$): the equivalence between $E(t)$ and $V(t)$ ; the estimation from below of $\dot V(t) = \frac {dV}{dt}(t)$ under the smart choice of the parameters $\gamma_0$, $\gamma_1$ of the event-triggering mechanism and of the weighting parameter $\eps$ ; the deduction of the  exponential decreasing of the energy $E$ of system \eqref{Eqet} under the smart choice of function $\eta_0(t)$ in \eqref{et-law} ; the conclusion in terms of the initial context.

\begin{proof} 

\noindent $\bullet$ {\it First Step:} The energy of the system and the proposed Lyapunov functional are equivalent if there exist two positive constants $C_1$ and $C_2$ such that for all $t\geq0$,
\begin{equation}\label{energy2}
C_1 E(t) \leq V(t) \leq C_2E(t).
\end{equation}
It is indeed easy to prove, using Cauchy-Schwartz, Young's and Poincar\'e 's inequalities recalled in appendix (with $\epsilon = C_\Omega$ for Lemma~\ref{Y}), that 
$$\begin{aligned}
&V(t) \leq E(t) + \dfrac {\eps\alpha}2  \|z(t)\|_{L^2(\Omega)}^2+ \eps \|z(t)\|_{L^2(\Omega)}\|\partial_t z(t)\|_{L^2(\Omega)}\\
&\leq  E(t) + \dfrac{\eps C_\Omega} 2 \|\partial_t z(t)\|_{L^2(\Omega)}^2 + \dfrac\eps 2 \left(\dfrac 1{C_\Omega} + \alpha\right) \| z(t)\|_{L^2(\Omega)}^2\\
&\leq  E(t) + \dfrac{\eps C_\Omega} 2 \|\partial_t z(t)\|_{L^2(\Omega)}^2 + \dfrac{\eps C_\Omega }{2 } \left( 1+ \alpha C_\Omega\right)\|\nabla  z(t)\|_{L^2(\Omega)}^2\\
&\leq  \left( 1+ \eps C_\Omega + \eps \alpha C_\Omega^2 \right) E(t),
\end{aligned}$$
so that $C_2 = ( 1+ \eps C_\Omega + \eps \alpha C_\Omega^2)$ in \eqref{energy2}.
On the other hand, for the same reasons, 
$$
V(t) \geq  E(t)  - \eps \|z(t)\|_{L^2(\Omega)}\|\partial_t z(t)\|_{L^2(\Omega)}\geq  \left( 1 - \eps C_\Omega \right) E(t),
$$
so that as soon as $\eps <  1/C_\Omega$, $C_1 = ( 1- \eps C_\Omega)$ in \eqref{energy2}.

\noindent $\bullet$ {\it Second Step:} The goal of this step is to prove that there exists $\beta>0$ such that for all $t\geq0$, 
\begin{equation}\label{Lyapunovderiv}
\dot V(t) \leq -\beta E(t) + \dfrac \alpha 2 (1+\alpha\eps) \eta_0(t).
\end{equation}
First, the closed-loop system \eqref{Eqet} with the event-triggering mechanism \eqref{et-law}  reads
\begin{equation}\label{event-triggered-wave}
	\left\{
	\begin{aligned}
		&\partial_t^2z(t,x) - \Delta z(t,x) = -\alpha \partial_t z(t,x)+\alpha e_k(t,x),\\
		& \qquad  \qquad  \qquad \qquad  ~  (t,x)\in[t_k,t_{k+1}[\times \Omega, \forall k\in\mathbb N\\
		&z(t,x)=0,\qquad \qquad \qquad \qquad (t,x)\in\mathbb{R}_+\times \partial \Omega\\
		&z(0,x)=z_0(x),\: \partial_t z(0,x)=z_1(x),\qquad \quad x\in \Omega.
	\end{aligned}
	\right.
\end{equation}
Compute the time-derivative of $V$ along the trajectories of system \eqref{event-triggered-wave} and use integrations by parts 
for all $t\geq0$, such that $t \in [t_k,t_{k+1}[$:
\begin{eqnarray*}
&&\dot V(t) =  \dfrac d{dt} \left( E(t) + \dfrac {\eps\alpha}2  \|z(t)\|_{L^2(\Omega)}^2+ \eps\int_\Omega  z(t)\partial_t z(t) \right)\\
& = & \int_\Omega \partial_t^2 z(t)\partial_t z(t) + \int_\Omega \partial_t\nabla  z(t)\cdot \nabla  z(t) \\
&&+~\eps\alpha \int_\Omega  z(t)\partial_t z(t)+ \eps\int_\Omega  |\partial_t z(t)|^2 + \eps\int_\Omega  z(t)\partial_t^2 z(t)
\end{eqnarray*}
\begin{eqnarray*}
& = & \int_\Omega \Delta z(t)\partial_t z(t) -\alpha \int_\Omega |\partial_t z(t)|^2 + \alpha \int_\Omega e_k(t)\partial_t z(t) \\
&&+ \int_\Omega \partial_t\nabla  z(t)\cdot \nabla  z(t)
+~ \eps\int_\Omega  |\partial_t z(t)|^2 + \eps \int_\Omega  z(t)\Delta z(t)\\
&&+~ \eps \alpha \int_\Omega  z(t)e_k(t)\\
& = & (\eps-\alpha) \int_\Omega |\partial_t z(t)|^2  - \eps\int_\Omega  |\nabla z(t)|^2 \\
&& + ~\alpha \int_\Omega e_k(t)\partial_t z(t) + \eps \alpha \int_\Omega  z(t)e_k(t).\\
\end{eqnarray*}
Now, let us estimate the two last terms of the equality using the definition of $e_k$ in \eqref{deviation} and the event-triggering mechanism \eqref{et-law} that gives, $\forall t \in [t_k,t_{k+1}[$
$$
\|e_k(t)\|_{L^2(\Omega)}^2 \leq \gamma_0  \|z(t)\|_{L^2(\Omega)}^2 + \gamma_1 \|\partial_t z(t)\|_{L^2(\Omega)}^2+ \eta_0(t).
$$
We obtain, using also Cauchy-Schwartz, Young and Poincar\'e's inequalities: 
$$\begin{aligned}
	\Big|&\alpha \int_\Omega e_k(t)\partial_t z(t)  \Big| 
	\leq \dfrac\alpha 2 \|e_k(t)\|_{L^2(\Omega)}^2  + \dfrac\alpha 2 \| \partial_t z(t)\|_{L^2(\Omega)}^2 \\
	&\leq  \dfrac{\alpha\gamma_0} 2   \|z(t)\|_{L^2(\Omega)}^2 +  \dfrac{\alpha(\gamma_1+1)} 2   \|\partial_t z(t)\|_{L^2(\Omega)}^2+ \dfrac\alpha 2 \eta_0(t)\\
	&\leq  \dfrac{\alpha\gamma_0C_\Omega^2} 2   \|\nabla z(t)\|_{L^2(\Omega)}^2 +  \dfrac{\alpha(\gamma_1+1)} 2   \|\partial_t z(t)\|_{L^2(\Omega)}^2\\
	&\hspace{7cm}+ \dfrac\alpha 2 \eta_0(t)
\end{aligned}$$
and
$$\begin{aligned}
	&\Big| \eps \alpha \int_\Omega  z(t)e_k(t) \Big| 
	\leq \dfrac{\eps}2 \| z(t)\|_{L^2(\Omega)}^2 + \dfrac{\alpha^2\eps}2 \| e_k(t)\|_{L^2(\Omega)}^2\\
	&\leq \dfrac{\eps C_\Omega^2}2 (1 + \alpha^2\gamma_0) \|\nabla z(t)\|_{L^2(\Omega)}^2 + \dfrac{\alpha^2\eps \gamma_1}2 \|\partial_t z(t)\|_{L^2(\Omega)}^2\\
	&\hspace{7cm}+ \dfrac{\alpha^2\eps}2\eta_0(t).
\end{aligned}$$
Therefore, gathering these estimations, we can write
$$
	\dot V(t) \leq - \dfrac {\nu_0}2 \|\nabla z(t)\|_{L^2(\Omega)}^2 - \dfrac {\nu_1}2 \|\partial_t z(t)\|_{L^2(\Omega)}^2  + \dfrac\alpha 2 \left( 1+ \alpha\eps \right) \eta_0(t)
$$
with
\begin{equation}\label{nu}
\left\{\begin{aligned}
&\nu_0 = 2 \eps  - \alpha\gamma_0C_\Omega^2 - \eps C_\Omega^2 (1 + \alpha^2\gamma_0)\\
&\nu_1 =  2 \alpha -2\eps - \alpha(\gamma_1+1)  -  \alpha^2\eps \gamma_1
\end{aligned}\right.
\end{equation}
and we aim at defining $\gamma_0$ and $\gamma_1$ such that we can choose $\eps>0$ to ensure $\nu_0 > 0$ and $\nu_1 > 0$.\\
\noindent $(i)$ To have $\nu_0>0$, we need 
$$
\eps \left( 2  -   C_\Omega^2 (1 + \alpha^2\gamma_0) \right) >  \alpha\gamma_0C_\Omega^2
$$
so that we have to define $\gamma_0 >0$ such that
\begin{equation}\label{gam}
 2  -   C_\Omega^2 (1 + \alpha^2\gamma_0) > 0,
\end{equation}
and then, we can choose $\eps$ satisfying
$
\eps > \dfrac{ \alpha\gamma_0C_\Omega^2}{2 -  C_\Omega^2 (1 + \alpha^2\gamma_0)}.
$
\\
Since we also have to satisfy $\eps < 1/C_\Omega$ from the end of the first step, we also need to check if 
\begin{equation}\label{gamm}
\dfrac{ \alpha\gamma_0C_\Omega^2}{2 -  C_\Omega^2 (1 + \alpha^2\gamma_0)} < \dfrac 1{C_\Omega}.
\end{equation}
Therefore, since \eqref{gam} is equivalent to  
$\gamma_0 <  \dfrac{2-  C_\Omega^2}{\alpha^2 C_\Omega^2}$
and \eqref{gamm} is equivalent to  
$
\gamma_0 <  \dfrac{2-  C_\Omega^2}{\alpha^2 C_\Omega^2 + \alpha C_\Omega^3},
$
the assumption $C_\Omega < \sqrt 2$ means that it is enough to set $\gamma_0$ satisfying 
\begin{equation}\label{gamma0}
0 < \gamma_0 <  \dfrac{2-  C_\Omega^2}{\alpha^2 C_\Omega^2+ \alpha C_\Omega^3} .
\end{equation}

\noindent $(ii)$ To have $\nu_1>0$, we need 
$
\eps \left( 2 + \alpha^2  \gamma_1\right) <  2 \alpha - \alpha(\gamma_1+1)
$
so that we have to define $\gamma_1< 1$ and then choose $\eps>0$ such that
$
\eps <  (\alpha (1-\gamma_1))/(2 + \alpha^2  \gamma_1).
$

In order to end this step, we have to assess the feasibility of the following estimate
\begin{equation}\label{Epsilon}
\dfrac{ \alpha\gamma_0C_\Omega^2}{2 -  C_\Omega^2 (1 + \alpha^2\gamma_0)} < \eps < \dfrac{\alpha (1-\gamma_1)}{2 + \alpha^2  \gamma_1}.
\end{equation}
Since $\gamma_0$ and $\gamma_1$ are such that both sides are positive, we only need to check that
$$
\alpha\gamma_0C_\Omega^2(2 + \alpha^2  \gamma_1) < \alpha (1-\gamma_1)(2 -  C_\Omega^2 (1 + \alpha^2\gamma_0))
$$
which is equivalent to 
$$
(\alpha^2-2)C_\Omega^2\gamma_0+   2 -  C_\Omega^2  + (C_\Omega^2-2) \gamma_1 > 0
$$
If $\alpha^2\geq  2$ then any $\gamma_0 >0$ and $\gamma_1<1$ are appropriate. If $\alpha^2 < 2$, then
it is easy to verify that if one defines actually 
\begin{equation}\label{gamma01}
\gamma_0 < \dfrac{2-C_\Omega^2}{2C_\Omega^2(2- \alpha^2)}~\hbox{ and }
\gamma_1 < \dfrac 12,
\end{equation}
estimate \eqref{Epsilon} can be satisfied.

Therefore, under appropriate conditions \eqref{gamma0} and \eqref{gamma01} to satisfy for $\gamma_0$ and $\gamma_1$, there exists an $\eps$ satisfying \eqref{Epsilon}, so that one can define $\beta>0$ from \eqref{nu} by 
$$
\beta = \min(\nu_0, \nu_1) 
$$
in order to get \eqref{Lyapunovderiv}.

\noindent $\bullet$ {\it Third Step:} From estimates \eqref{energy2}  and \eqref{Lyapunovderiv}, one can write that for all $t\in[t_k,t_{k+1}[$,
$$
\dot V(t) \leq -\dfrac \beta{C_2} V(t) + \dfrac \alpha 2 (1+\eps) \eta_0(t),
$$
which gives 
$
\dfrac{d}{dt} \left( V(t) e^{\frac \beta{C_2}t }\right) \leq \dfrac \alpha 2 (1+\eps) \eta_0(t)e^{\frac \beta{C_2}t }.
$\\
Let us now choose
\begin{equation}\label{eta}
\eta_0(t) = V(0)e^{-\theta t} \hbox{ with } \theta > \dfrac{\beta}{1+\eps C_\Omega+\eps \alpha C_\Omega^2}.
\end{equation}
Then, for all $t\in]t_k,t_{k+1}]$, by integration between $t_k$ and $t$ one obtains
$$
 V(t) e^{\frac \beta{C_2}t }  -  V(t_k) e^{\frac \beta{C_2}t_k }  \leq \dfrac \alpha 2 (1+\eps)  V(0) \int_{t_k}^t e^{-\theta s}e^{\frac \beta{C_2}s} ds.
$$
Now, we calculate easily for all $t\in]t_k,t_{k+1}]$ 
$$ V(t) \leq  \mu  V(0) \left( e^{-\theta t_k}e^{-\frac \beta{C_2}(t-t_k)} - e^{-\theta t} \right)  
  + V(t_k) e^{-\frac \beta{C_2}(t-t_k)}
$$
where $\mu = \dfrac {\alpha (1+\eps)}{2(\theta - \frac \beta{C_2})} >0$ since $\theta > \dfrac \beta{C_2}$.\\
Therefore, reasoning by recurrence on $k$ (down to $t_0 = 0$), we can derive that for any $t>0$, there exists a unique $k\in \mathbb N$ such that $t \in [t_k,t_{k+1}[$ and 
\begin{eqnarray*}
 V(t) &\leq&   \mu  V(0) \left( e^{-\theta t_k}e^{-\frac \beta{C_2}(t-t_k)} - e^{-\theta t} \right)  \\
 && +~ \mu  V(0) \left( e^{-\theta t_{k-1}}e^{-\frac \beta{C_2}(t_k-t_{k-1})} - e^{-\theta t_k} \right) e^{-\frac \beta{C_2}(t-t_k)}\\
&&  + ~V(t_{k-1}) e^{-\frac \beta{C_2}(t_k-t_{k-1})}  e^{-\frac \beta{C_2}(t-t_k)}\\
 &\leq&    \mu  V(0) \left(  - e^{-\theta t} + e^{-\theta t_{k-1}}e^{-\frac \beta{C_2}(t-t_{k-1})}    \right)  \\
&&  +~ V(t_{k-1}) e^{-\frac \beta{C_2}(t-t_{k-1})}  \\
&\leq&   (1+\mu ) V(0) e^{-\frac \beta{C_2}t}   - \mu  V(0)e^{-\theta t}.
\end{eqnarray*}
Recalling \eqref{energy}, in terms of the energy $E$ of the system, we proved
$$ E(t) \leq  \dfrac {C_2}{C_1}  (1+\mu ) E(0) e^{-\frac \beta{C_2}t}   - \mu  E(0)e^{-\theta t}$$
which gives finally
$ E(t) \leq K E(0) e^{-\delta t}
$
with the overshoot estimation $K$:
\begin{equation}\label{KK}
K =   \dfrac {1+\eps C_\Omega+\eps\alpha C_\Omega^2}{1-\eps C_\Omega}  \left[1+  \dfrac {\alpha (1+\eps)}{2(\theta - \frac \beta{1+\eps C_\Omega+\eps\alpha C_\Omega^2})} \right] 
\end{equation}
and the exponential rate $\delta$:
\begin{equation}\label{dd}
\delta = \frac \beta{1+\eps C_\Omega+\eps\alpha C_\Omega^2}.
\end{equation}

\noindent $\bullet$ {\it Fourth Step:} 
The goal of this step is to state the conclusion of the proof of Theorem~\ref{GES} in terms of the initial context of system \eqref{Eqet} and event-triggering mechanism \eqref{et-law}, based only on  $\gamma_0$, $\gamma_1$, $C_\Omega$ and $\alpha$.

Gathering all the information about the appropriate choice of the parameters introduced in the proof, namely $\eps$, $\beta$ and $\theta$, it follows: $\eps$ has to satisfy 
$\eps <  1/C_\Omega$ and \eqref{Epsilon} ;  $\beta = \min(\nu_0, \nu_1)$ with \eqref{nu} and finally $\theta >  \beta/C_2$ as in \eqref{eta}.

The possibility of choosing the parameters ($\gamma_0$, $\gamma_1$ and $\eta_0$) of the event-triggering mechanism \eqref{et-law}  results from the conditions \eqref{gamma0}, \eqref{gamma01} and \eqref{eta}. It can be summarized as follows: 
$$
\left\{ \begin{array}{lcl}
\begin{array}{c}
\gamma_0 < \min \left(   \dfrac{2-  C_\Omega^2}{\alpha^2 C_\Omega^2+ \alpha C_\Omega^3},
\dfrac{2-C_\Omega^2}{2C_\Omega^2(2- \alpha^2)}\right) \\
\qquad \qquad \qquad \hfill{} \mbox{ and } 
\gamma_1 <  1/2 \end{array} 
& \mbox{if} & \alpha < 2\\
\gamma_0 < \dfrac{2-  C_\Omega^2}{\alpha^2 C_\Omega^2+ \alpha C_\Omega^3}
\mbox{ and } \gamma_1 < 1 & \mbox{if} & \alpha \geq 2
\end{array}\right.
$$
 and
$$
 \eta_0(t) = \left( \frac{1}{2}\|z_1\|_{L^2(\Omega)}^2 + \frac{1}{2} \|\nabla  z_0\|_{L^2(\Omega)}^2 
+ \eps\int_\Omega  z_0z_1 \right) e^{-\theta t } 
$$
with $\theta > \dfrac{\min(\nu_0, \nu_1)}{1+\eps C_\Omega+\eps \alpha C_\Omega^2}$,
concluding Theorem~\ref{GES}. 
\end{proof}

From the proof of Theorem \ref{GES}, one can state the following result, presenting explicit numerical
conditions to design the adequate even-triggering mechanism  \eqref{et-law}.
\begin{theorem}\label{GES2}
Given  the damping parameter $\alpha>0$ in system \eqref{Eqet} and 
given the Poincar\'e constant $C_\Omega$ satisfying $C_\Omega < \sqrt 2$, if  
$\gamma_0>0$, $\gamma_1>0$ and $\eta_0 \in C_0(\mathbb R_+)$ satisfy:
$$
\gamma_0 <  \dfrac{2-  C_\Omega^2}{ C_\Omega^2 \max(\alpha^2 + \alpha C_\Omega, 2(2- \alpha^2))}
\mbox{ and } 
\gamma_1 < \dfrac 12 \mbox{, if } \alpha < 2$$
or 
$$
\gamma_0 < \dfrac{2-  C_\Omega^2}{\alpha^2 C_\Omega^2+ \alpha C_\Omega^3}
\mbox{ and }  \gamma_1 < 1 \mbox{, if } \alpha \geq 2,
$$ 
 and also
$$
 \eta_0(t) = \left( \frac{1}{2}\|z_1\|_{L^2(\Omega)}^2 + \frac{1}{2} \|\nabla  z_0\|_{L^2(\Omega)}^2 
+ \eps\int_\Omega  z_0z_1 \right) e^{-\theta t } 
$$
with $\theta > \min(\nu_0, \nu_1) / \left(1+\eps C_\Omega+\eps \alpha C_\Omega^2\right)$, \\
then
for any initial condition $(z_0,z_1)\in H^1_0(\Omega)\times L^2(\Omega)$, system \eqref{Eqet} under the event-triggering mechanism \eqref{et-law} is exponentially stable with the overshoot estimation $K$ and the exponential rate $\delta$ given by \eqref{KK} and \eqref{dd}
where $\eps$ satisfies $\eps <  1/C_\Omega$ and \eqref{Epsilon}, and $\beta$ is defined from \eqref{nu} by $\beta = \min(\nu_0, \nu_1)$.

\end{theorem}

From this theorem, it should be interesting to be able to orient the choice of the parameters $\gamma_0$ and $\gamma_1$, thanks to an optimization scheme. This constitutes an ongoing study with the goal to reduce the number of control updates.
\section{Discussion and conclusion}
\label{sec_discussion}
We have designed an event-triggered control for a linear wave equation. After proving that the associated event-triggering mechanism  avoids Zeno solutions, we provided a Lyapunov stability analysis proving that the system is exponentially stable. Explicit expressions allowing to choose the parameters characterizing the event-triggering rule have been proposed.

These results pave the way to many other works to be done. Let us mention some of them:\\
$\bullet$ One first step would consist in avoiding the restriction on the size of the domain $\Omega$ (see Remark \ref{rem-restriction-domain} for more details). This could be possible for instance by focusing on other methods than Lyapunov ones, for instance the multiplier method as in \cite{komornik1997exact}, allowing to find sharp constants for the stabilization of PDEs.\\
$\bullet$ Studying the case of a localized damping is challenging, since it is often impossible to provide a Lyapunov functional in this situation. One uses rather a compactness/uniqueness strategy, that has been applied either for the wave equation (\cite{saut1987unique}) or the Korteweg-de Vries equation (\cite{rosier2006global}). Such a strategy could be followed for the sampled localized or boundary damping of a linear wave equation.

\section{Appendix}
\begin{lemma}[Cauchy-Schwartz inequality]
$$
\int_\Omega u(x)v(x)dx \leq \|u\|_{L^2(\Omega)}\|v\|_{L^2(\Omega)}, \qquad
\forall u, v \in L^2(\Omega).
$$
\end{lemma} 
\begin{lemma}[Young's inequality]\label{Y}
$$
ab \leq \dfrac{\epsilon}2 a^2 + \dfrac{1}{2\epsilon} b^2, \qquad \forall a,b \in \mathbb R, \forall \epsilon>0.
$$
\end{lemma} 
\begin{lemma}[Poincar\'e's inequality]\label{Poincare}
There exists a constant $C_\Omega>0$ depending on the size, the geometry and the regularity of the bounded domain $\Omega$ such that for any function $u \in H^1_0(\Omega)$, 
$
\|u\|_{L^2(\Omega)} \leq C_\Omega \|\nabla u\|_{L^2(\Omega)}.
$
\end{lemma} 
The optimal constant $C_\Omega$ in the Poincar\'e inequality is sometimes known as the Poincar\'e constant for the domain $\Omega$.
For instance, if the space dimension is $N=1$, Wirtinger's inequality gives $C_{[a,b]} =  (b-a)/(2\pi)$ and in general, 
if $\Omega \subset \mathbb R^N$ is a convex Lipschitzian domain of diameter $d$, then $C_\Omega \leq  d/\pi$. 


\end{document}